\documentclass{gtart}     
\newtheorem{thm}{Theorem}[section]

\newtheorem{lem}[thm]{Lemma}
\newtheorem{prop}[thm]{Proposition}
\def\bdy{\partial}
\def\rel{\hbox{$\,\hbox{\it rel\/}\,$}}
\def\diff{\hbox{\it diff\/}}
\def\Diff{\hbox{\it Diff\/}}
\def\BDiff{\hbox{\it BDiff\/}}
\def\hscr{{\cal H}}

\def\Imb{\hbox{\it Imb\/}}
\def\Isom{\hbox{\it Isom\/}}
\def\Out{\hbox{\it Out\/}}
\def\Z{\hbox{\Bbb Z}}

\font\Bbb=msbm10   scaled 1096
\def\proclaim #1. #2\par{\ppar{\bf#1}\stdspace{\sl#2}\ppar}
\begin{document}

%
%
\input gtoutput
\volumenumber{1}\papernumber{7}\volumeyear{1997}
\pagenumbers{91}{109}\published{21 December 1997}
\shorttitle{Diffeomorphism Groups of 3-Manifolds}  
\proposed{Robion Kirby}\seconded{Joan Birman, David Gabai}
\received{12 June 1997}\revised{19 December 1997}

\title{Finiteness of Classifying Spaces of\\Relative Diffeomorphism Groups 
of 3--Manifolds}
\shorttitle{Classifying spaces of difeeomorphism groups of 3--manifolds}

\authors{Allen Hatcher\\Darryl McCullough}  

\address{Department of Mathematics, Cornell University\\
Ithaca, NY 14853, USA}           
\secondaddress{Department of Mathematics, University of Oklahoma\\
Norman, OK 73019, USA}
\email{hatcher@math.cornell.edu}
\secondemail{dmccullough@math.ou.edu}

\asciiaddress{Department of Mathematics, Cornell University,
Ithaca, NY 14853, USA,
Department of Mathematics, University of Oklahoma
Norman, OK 73019, USA,
hatcher@math.cornell.edu, dmccullough@math.ou.edu}

\begin{abstract}

The main theorem shows that if $M$ is an irreducible compact
connected orientable 3--manifold with non-empty boundary, then the
classifying space\break $\BDiff(M\rel \partial M)$ of the space of
diffeomorphisms of $M$ which restrict to the identity map on $\partial
M$ has the homotopy type of a finite aspherical CW--complex. This
answers, for this class of manifolds, a question posed by
M~Kontsevich. The main theorem follows from a more precise result,
which asserts that for these manifolds the mapping class group
$\hscr(M\rel \partial M)$ is built up as a sequence of extensions of
free abelian groups and subgroups of finite index in relative mapping
class groups of compact connected surfaces.

\end{abstract}

\asciiabstract{The main theorem shows that if M is an irreducible
compact connected orientable 3-manifold with non-empty boundary, then
the classifying space BDiff(M rel dM) of the space of diffeomorphisms
of M which restrict to the identity map on dM has the homotopy type of
a finite aspherical CW-complex. This answers, for this class of
manifolds, a question posed by M Kontsevich. The main theorem follows
from a more precise result, which asserts that for these manifolds the
mapping class group H(M rel dM) is built up as a
sequence of extensions of free abelian groups and subgroups of finite
index in relative mapping class groups of compact connected surfaces.}





\maketitlepage

%

For a compact connected 3--manifold $M$, let $\Diff(M\rel R)$ denote
the group of diffeomorphisms $M \to M$ restricting to the identity on
the subset $R$. We give $\Diff(M\rel R)$ the $C^\infty$--topology, as
usual. M.~Kontsevich has conjectured (problem~3.48
in~\cite{Kirbylist}) that the classifying space $\BDiff(M\rel\partial
M)$ has the homotopy type of a finite complex when $\bdy M$ is
non-empty.  In this paper we prove the conjecture for irreducible
orientable 3--manifolds. In fact, more is true in this case.

\proclaim Main Theorem. Let $M$ be an irreducible compact connected
orientable 3--manifold and let $R$ be a non-empty union of components of
$\bdy M$, including all the compressible ones. Then $\BDiff(M\rel R)$
has the homotopy type of an aspherical finite CW--complex.

\noindent
Actually, the assertion that $\BDiff(M\rel R)$ is aspherical has been
known for some time (see \cite{H1,H2,I}), as a special case of more
general results about Haken manifolds. Thus the finiteness question
for $\BDiff(M\rel R)$ is equivalent to whether the mapping class group
$\pi_0(\Diff(M\rel R))$, which we denote by\break $\hscr(M\rel R)$, is a
group whose classifying space is homotopy equivalent to a finite
complex. Such groups are called {\it geometrically finite.}

It is a standard elementary fact that a geometrically finite group
must be torsion-free. Thus the Main Theorem implies that $\hscr(M\rel
R)$ is torsion-free, a fact which can be deduced
from~\cite{Heil-Tollefson}. If we drop the condition that the
diffeomorphisms restrict to the identity on $R$, or allow $M$ to be
closed, then $\hscr(M)$ can have torsion (for example if $M$ is a
hyperbolic 3--manifold with non-trivial isometries), and thus
$\BDiff(M)$ can be aspherical but not of the homotopy type of a finite
complex. For all Haken 3--manifolds, however, there exist
finite-sheeted covering spaces of $\BDiff(M)$ and $\BDiff(M\rel R)$
which have the homotopy type of a finite complex~\cite{M}. The Main
Theorem can be viewed as a refinement of this result for the case
of~$\BDiff(M\rel R)$.

If the irreducibility condition on $M$ is dropped, $\BDiff(M\rel R)$
need no longer be aspherical. Indeed, its higher homotopy groups can
be rather complicated, and in particular $\pi_2(\BDiff(M\rel R))$ is
generally not finitely generated~\cite{K-M}. This does not exclude the
possibility that $\BDiff(M\rel R)$ has the homotopy type of a finite
complex (for example, $S^1 \vee S^2$ is a finite complex having
non-finitely generated~$\pi_2$), and it would be very interesting to
know whether Kontsevich's conjecture holds in these cases. If so, it
would indicate that $\BDiff(M\rel R)$ is more tractable than has
generally been supposed.

The Main Theorem follows directly from a structural result about
$\hscr(M\rel R)$ which we will state below. To place our result in
historical context, and to review some of the ingredients that go into
its proof, we will first survey some previous work on mapping class
groups of Haken 3--manifolds. Recall that according to the basic
structure theorem of Jaco--Shalen~\cite{JS} and Johannson~\cite{J2}, a
Haken manifold with incompressible boundary admits a characteristic
decomposition into pieces which are either fibered (admitting an
$I$--fibering or a Seifert fibering) or simple (admitting no
incompressible torus or proper annulus which is not properly homotopic
into the boundary). By the late 1970's, substantial information had
been obtained about the mapping class groups of these pieces. For the
pieces that are $I$--bundles, the mapping class group is essentially
the same as the mapping class group of the quotient surface. The
necessary technical results to analyze this case are contained in
Waldhausen's seminal paper~\cite{W}. For a Seifert-fibered piece $W$,
the analysis was due to Waldhausen (pages~85--86 of~\cite{W}, page~36
of~\cite{PSPM}) and Johannson (proposition 25.3 of~\cite{J2}). After
showing that one can restrict attention to diffeomorphisms preserving
the fiber structure, they deduced that $ \hscr(W) $ fits into a short
exact sequence where the kernel group is finitely-generated abelian
and the quotient group is a surface mapping class group. For the
simple pieces, Johannson proved that the group of mapping classes
preserving the frontier is finite (proposition~27.1 of~\cite{J2}). Of
course, this was carried out without reference to the hyperbolic
structure later discovered to exist on these pieces. Today, this
finiteness is often viewed as a consequence of Mostow rigidity, which
implies that if $W$ is a 3--manifold with a complete hyperbolic
structure with finite volume, then $\hbox{\sl Out}(\pi_1(W))$ is
finite, and from Waldhausen's fundamental work, this group is
isomorphic to the mapping class group. However, when $\partial M$ has
components other than tori, the simple pieces of $M$ might not admit
hyperbolic structures of finite volume, and indeed their mapping class
groups may not be finite, but their groups of mapping classes
preserving the frontier will be finite.

By combining the information on these two types of pieces, Johannson
proved the first general result on mapping class groups of Haken
3--manifolds. This result, corollary~27.6 in~\cite{J2}, says that in
the case when $M$ has incompressible boundary, the subgroup of mapping
classes generated by Dehn twists about tori and properly imbedded
annuli in $M$ has finite index in~${\cal H}(M)$. (For a definition of
Dehn twist, see section~3 below.)

At about the same time, techniques for controlling isotopies between
diffeomorphisms of 3--manifolds were being developed, leading to more
refined structural details about mapping class groups.
Laudenbach~\cite{Lau} proved that (apart from a few easily-understood
exceptions) an isotopy between two diffeomorphisms of a Haken
3--manifold that preserve an incompressible surface can be deformed to
an isotopy that preserves the surface at each level of the isotopy.
This was extended to parameterized families by Hatcher~\cite{H2} and
Ivanov~\cite{I}, giving the key ingredient in the proof that when
$\partial M$ is non-empty, the components of $\Diff(M\rel\partial M)$
are contractible, and consequently $\BDiff(M\rel\partial M)$ is an
aspherical complex. Laudenbach's result also led to the extension of
Harer's homological finiteness results on 2--manifold mapping class
groups to dimension~3 by McCullough~\cite{M}.

Among Harer's results was the fact that 2--manifold mapping class
groups contain geometrically finite subgroups of finite index. Using
Harer's constructions, in the simplified exposition of~\cite{H3}, we
strengthen this as follows.

\proclaim Lemma~\ref{HARER}. Let $S$ be a compact connected surface,
and let $J$ and $K$ be 1--dimensional submanifolds of $\partial S$ with
$J\cap K=\partial J\cap\partial K$. If $K$ is non-empty, then
$\hscr(S,J\rel K)$ is geometrically finite.

\noindent Here, the notation $\hscr(S,J\rel K)$ indicates the mapping
classes that carry $J$ diffeomorphically to $J$. The geometric
finiteness of $\hscr(M\rel R)$ is then an immediate consequence of the
following structure theorem, whose proof occupies most of this paper.

\proclaim Filtration Theorem.
Let $M$ be an irreducible compact connected orientable 3--man\-i\-fold
and let $R$ be a non-empty union of components of $\bdy M$,
including all the compressible ones. Then there is a filtration
$$0 = G_0 \subset G_1 \subset \cdots \subset G_n =
\hscr (M\rel R)
$$
\noindent
with each $G_i$ a normal subgroup of $G_{i+1}$ such that $G_{i+1}/G_i$
is either a finitely generated free abelian group or a subgroup of
finite index in a mapping class group $\hscr(S,J\rel K)$ of a compact
connected surface~$S$.

Let us describe how the filtration arises. An argument using a family
of compressing discs for the compressible components of $\partial M$
reduces the proof to the case that $\partial M$ is incompressible.
Then, each mapping class contains representatives which preserve the
characteristic decomposition of Jaco--Shalen and Johannson, and
Laudenbach's result implies that isotopies between such
diffeomorphisms can also be assumed to preserve the decomposition. We
can filter $M$ by an increasing sequence of submanifolds $V_i$, each
obtained from the preceding one by attaching one of the characteristic
pieces. This leads to a decreasing filtration of $\hscr(M \rel R)$ by
the subgroups represented by diffeomorphisms which restrict to the
identity on succesively larger $V_i$'s. The preceding remarks imply
that the successive quotients for this filtration are subgroups of
relative mapping class groups of the corresponding attached
characteristic pieces. This reduces us to understanding the relative
mapping class groups of the characteristic pieces.

In the case of fibered pieces, the analysis is that of Waldhausen and
Johannson. We present it in the relative cases that we will need as
lemmas~\ref{ibundle mcg} and~\ref{seifert mcg}. The upshot is that
after interpolating one more stage in the filtration, the two
resulting quotient groups are of the desired types. For the simple
pieces, the finiteness result of Johannson implies that the relative
mapping class group has a finitely-generated abelian subgroup of
finite index, generated by Dehn twists along peripheral tori and
annuli. But to rule out the existence of torsion, we show that the
relative mapping class group itself is finitely-generated free
abelian. This is done using an appropriate form of Mostow rigidity. In
the case that $\bdy M$ consists entirely of tori, the usual form of
Mostow rigidity suffices, but in the general case an extended version
is needed. To state this, let $W$ be a compact connected orientable
irreducible 3--manifold with non-empty boundary, let $T$ be the union of
its torus boundary components, and let $A$ be a union of disjoint
incompressible annuli in $\partial W-T$. Certain assumptions, listed
at the beginning of section~3, are made which are satisfied when $W$
is a simple piece of the characteristic decomposition of a Haken
3--manifold and $A$ is the union of the components of the frontier of
$W$ that are annuli. Let $\Diff(W,A)$ denote the diffeomorphisms of
$W$ that take $A$ diffeomorphically to $A$, and let ${\cal
H}(W,A)=\pi_0(\Diff(W,A))$.

\proclaim Lemma~\ref{HYPERBOLIC STRUCTURE}. $W-(A\cup T)$ has a
hyperbolic structure with totally geodesic boundary. Its group of
isometries $\Isom(W-(A\cup T))$ is finite, and $\Isom(W-(A\cup
T))\to\hscr(W,A)$ is an isomorphism.

\noindent This version of Mostow Rigidity is folklore, but as will be
seen, it is not such a simple matter to give a real proof. A key
ingredient is a theorem of Tollefson~\cite{Tollefson}, which provides
a very strong uniqueness statement for certain involutions of Haken
3--manifolds.

\section{Preliminaries}
\label{preliminaries}

In this section we collect some auxiliary results that will be used in
the proof of the Main Theorem.

As explained in the introductory section, the Filtration Theorem will
show that $\hscr(M\rel R)$ is built up by a sequence of extensions of
geometrically finite groups. Then, the following lemma will imply the
Main Theorem.

\begin{lem}\sl\label{extensions of gf groups} Let $1 \to H \to G \to K
\to 1$ be a short exact sequence of groups. If $H$ and $K$ are
geometrically finite, then $G$ is geometrically finite.
\end{lem}

\begin{proof} By proposition 5(c) of \cite{Serre}, a group $\Gamma$
is geometrically finite if and only if it is finitely presented and
FL. The latter means that there is a finite length resolution of the
trivial $\Z\Gamma$--module $\Z$ by finitely generated free
$\Z\Gamma$--modules. By proposition~6(b) of~\cite{Serre}, extensions of
FL groups are FL, and the lemma follows.
\end{proof}

The next lemma provides geometrically finite groups that will form
some of the building blocks for $\BDiff(M\rel R)$.

\begin{lem}\sl\label{HARER} Let $S$ be a compact connected surface,
and let $J$ and $K$ be 1--dimensional submanifolds of $\partial S$ with
$J\cap K=\partial J\cap\partial K$. If $K$ is non-empty, then
$\hscr(S,J\rel K)$ is geometrically finite.
\end{lem}

\begin{proof}
We use constructions due to Harer, in the simplified exposition of
\cite{H3}. We may assume that $J$ consists of arcs, since if $J'$
consists of the arc components of $J$, then ${\cal H}(S,J\rel K)$ has
finite index in ${\cal H}(S,J'\rel K)$. As a further simplification,
we may assume all the components of $K$ are circles, since replacing
each arc of $K$ by the component of $\bdy S$ containing it, and
deleting from $J$ any components engulfed by $K$ in this process, has
no effect on ${\cal H}(S,J\rel K)$.

Let $P$ be the finite set obtained by choosing one point from each
circle of $K$ and one point from the interior of each arc of $J$.
Consider finite systems of arcs in $S$ with endpoints in $P$ and with
interiors disjointly embedded in the interior of $S$, such that:

\begin{enumerate}
\item[(a)]  Each arc is essential:  cutting $S$ along the arc does not
produce two components, one of which is a disc. This is equivalent to
saying that the arc represents a non-trivial element of  $\pi_1(S,\bdy
S)$.
\item[(b)] No two arcs in a system are isotopic in $S$ rel endpoints.
\end{enumerate}

\noindent
If $S$ is a disc, or annulus for which $J\cup K$ meets only one
component of $\partial S$, then $\hscr(S, J\rel K)$ is trivial and the
lemma holds. Otherwise, form a simplicial complex $\cal A$ whose
$k$--simplices are the isotopy classes of systems of $k+1$ arcs
satisfying (a) and (b). The barycentric subdivision $\cal A'$ of $\cal
A$ is the simplicial complex associated to the partially ordered set
of isotopy classes of arc systems, with the partial ordering given by
inclusion of systems. We are interested in the subcomplex $\cal B
\subset \cal A'$ associated to the partially ordered set of systems
whose complementary components are either discs or once-punctured
discs, the puncture being a component of $\bdy S$ that does not meet
$J\cup K$. The proof in~\cite{H3} describes a surgery process that
determines a flow on $\cal A$ that moves $\cal A$ into the star of a
vertex corresponding to a single arc. This flow preserves $\cal B$. If
the process is successively repeated for the arcs that make up a
vertex of $\cal B$, then $\cal B$ will be moved into the star of that
vertex, hence can be further contracted to the vertex moving only
through $\cal B$. Thus $\cal B$ is also contractible.

The group $\hscr (S,J \rel K)$ acts simplicially on $\cal A$ and $\cal
B$.  The action on $\cal B$ is without fixed points.  For if a point
in a simplex of $\cal B$ were fixed by an element of $\hscr (S,J \rel
K)$, the simplex would be invariant, hence fixed since its vertex arc
systems are distinguished from each other by the number of arcs they
contain. Thus one would have an arc system $ \alpha $, representing a
vertex of $\cal B$, which is taken by an element  $h \in \Diff
(S,J\rel K)$ to an isotopic arc system $h(\alpha)$. By isotopy
extension, we may assume $ h(\alpha)=\alpha$. The defining property of
${\cal B}$ implies that each component of $S-\alpha$ is either a disc
or a disc with a puncture corresponding to a component of $\partial S$
that does not meet $J\cup K$. Since $K$ is non-empty, $h$ must preserve
the closure of at least one of the components of $S-\alpha$ and fix
all elements of $P$ that it contains. By induction on the number of
components of $S-\alpha$, we deduce that $h$ must fix each point of
$P$ and preserve each component of $S-\alpha$. It follows that $h$ is
isotopic, relative to $K$ and preserving $J$, to the identity of~$S$.

Thus the quotient $ {\cal B} / \hscr (S,J \rel K)$ is a $K(\hscr (S,J
\rel K),1)$.  This quotient is a finite complex since arc systems fall
into finitely many orbits under the action of $\Diff (S,J \rel K)$.
\end{proof}

We will use the following consequence of a theorem of Laudenbach.

\begin{lem}\sl\label{preserve F} Let $M$ be a compact connected
irreducible 3--manifold which does not contain two-sided projective
planes. In $M$ let $F$ be a properly imbedded 2--sided incompressible
2--manifold, no component of which is a 2--sphere. Let $J_t\co M\to M$
be an isotopy such that $J_0$ is the identity and $J_1(F)= F$. If
either
\begin{enumerate}
\item[{\rm(i)}] $\partial F$ is non-empty and $J_t(\partial F)=\partial
F$ for all $t$, or
\item[{\rm(ii)}] $M$ does not fiber over $S^1$ with $F$ as fiber,
\end{enumerate}
\noindent then $J_t$ is deformable (through isotopies and relative to
$M\times\partial I$) to an isotopy which preserves $F$ at each level.
In case~(i), or in case~(ii) when $F$ is closed, the deformation can
also be taken relative to $\partial M\times I$.
\end{lem}

\begin{proof}
By induction we may assume that $F$ is connected. Fix a basepoint
$x_0$ in the interior of $F$. The proof of theorem~7.1 of~\cite{W}
shows that under hypothesis~(i) or~(ii), $J_t$ is deformable relative
to $M\times \partial I$ to a homotopy $H_t$ which preserves $F$ at
each level, and also shows that $J_1\vert_F$ must be an
orientation-preserving diffeomorphism of $F$. Let $h_t$ be the
restriction of $H_t$ to $F$. Any homotopy from the identity map to an
orientation-preserving diffeomorphism of $F$ can be deformed relative
to $F\times I$ to an isotopy~\cite{G}, that is, there exists an
isotopy $h'_t$ from the identity of $F$ to $h_1$, such that the path
$h_t$ followed by the reverse of $h'_t$ is a contractible loop in the
space of homotopy equivalences of $F$. Let $K_t$ be the isotopy of $M$
obtained by extension of the reverse of $h'_t$, starting from $J_1$,
and let $L$ be the product isotopy $JK$. Then $L$ has trivial trace at
$x_0$, and $L_1$ is the identity on $F$. Since $\Diff(M\rel\partial
M)\to\Imb(x_0,M-\partial M)$ is a fibration, $L$ is deformable
relative to $\partial M\times I$ to an isotopy that fixes~$x_0$.

Let $\ell_t$ be the restriction of $L_t$ to $F$.  What is proven on
pages~49--62 of~\cite{Lau} (see the comments at the end of page~48) is that
$\pi_1(\Imb(F,M\rel x_0))=0$. So $\ell_t$ is deformable to the
constant loop at the inclusion. Since $\Diff(M)\to\Imb(F,M)$ is a
fibration, this deformation of $\ell_t$ extends to a deformation of
$L_t$ to an isotopy which is the identity on $F$ for every $t$. Since
$K$ preserves $F$ at every level, it follows that $J$ is deformable to
an isotopy $J'$ which preserves $F$ at every level. When $F$ is
closed, all deformations can be taken relative to $\partial M\times
I$. In case~(i), the trace of $J$ at a point in $\partial F$ is a path
in $\partial F$. Let $G$ be a component of $\partial M$. Using
\cite{G}, and the fact that $G$ is not the 2--sphere, any two paths in
$\Diff(G)$ with the same trace and the same endpoints are deformable
to each other. Therefore $J'$ can be deformed to agree with $J$ on
$\partial M\times I$. Since $\pi_2\Diff(G)= 0$, the deformation from
$J$ to $J'$ can then be taken relative to $\partial M\times I$.
\end{proof}

\section{Fibered Manifolds}
\label{fibered}

In this section $M$ will be a compact connected orientable 3--manifold
whose boundary is decomposed as the union of two compact subsurfaces
$A$ and $B$ which intersect only in the circles of $\bdy A = \bdy B$.
We assume that the components of $A$ are annuli. We shall be
considering $\Diff(M,A\rel R)$, the group of diffeomorphisms of $M$
taking $A$ to itself and restricting to the identity on $R$, a
non-empty union of components of $A$ and $B$.

Suppose first that $M$ is an $I$--bundle over a compact connected
surface $S$, with projection map $p \co M \to S $.  We let $A$ be
the union of the fibers over $\bdy S$, so $B$ is the associated
$\partial I$--bundle.

\begin{lem}\sl\label{ibundle mcg} If $S$ has negative Euler
characteristic, then $\hscr(M,A\rel R) $ is isomorphic to $
\hscr(S\rel p(R))$. In particular, $\hscr(M,A\rel R) = 0$ if $R$ is
not contained in $A$.
\end{lem}

\begin{proof} We assert that the inclusion $\Diff_f(M,A\rel R)
\hookrightarrow \Diff(M,A\rel R)$ of the subgroup consisting of
diffeomorphisms taking fibers to fibers induces an isomorphism on
$\pi_0$. In the case when $R$ contains at most one component of $B$,
this follows from corollary~5.9 of~\cite{J2}. Otherwise, $M$ is a
product $I$--bundle and $B$ has two components. Let $R'=\overline{R-L}$
where $L$ is one of the components of $B$. Then the assertion holds
for $\Diff(M,A\rel R')$ and shows that it is contractible. Since $L$
has negative Euler characteristic, the identity component $\diff(L\rel
L\cap R')$ is contractible. Therefore the fibration $\Diff(M,A\rel
R')\to\diff(L\rel L\cap R')$ is a homotopy equivalence, so its fiber
$\Diff(M,A\rel R)$ is contractible and the assertion holds in this
case as well.

Now, viewing $M$ as the mapping cylinder of the projection $B \to S$,
there is a subgroup of $\Diff_f(M\rel R)$ consisting of
diffeomorphisms taking each level $B\times \{ t\} $ of the mapping
cylinder to itself, and the inclusion of this subgroup also induces an
isomorphism on $\pi_0$ since $\Diff(I\rel\bdy I)$ is contractible and
$R$ is non-empty. Diffeomorphisms in this subgroup are determined by
the quotient diffeomorphism they induce on $S$, and the result
follows.
\end{proof}

Suppose now that $M$ is an orientable compact connected irreducible
3--manifold Seifert fibered over the surface $S$, with projection
$p\co M\to S$. We assume the annuli of $A$ in the decomposition
$\bdy M = A\cup B$ are unions of fibers.

The images of the exceptional fibers form a finite set of exceptional
points $E \subset S - \bdy S$.  Each exceptional point can be labelled
by a rational number normalized to lie in the interval $(0,1)$,
describing the local structure of the Seifert fibering near the
corresponding exceptional fiber of $M$.  Let $\Diff^*(S,E\cup p(A)
\rel p(R))$ be the subgroup of $\Diff(S,p(A)\rel p(R))$ consisting of
diffeomorphisms permuting the points of $E$ in such a way as to
preserve the labelling, and let $\hscr^*(S,E\cup p(A) \rel p(R))$
denote the corresponding mapping class group.

\begin{lem}\sl\label{seifert mcg}  There is a split short
exact sequence
$$
0   \rightarrow   H_1 (S,\bdy S-p(R);{\Z})   \rightarrow
\hscr(M,A\rel R)    \rightarrow
\hscr^*(S,E\cup p(A) \rel p(R))
\rightarrow   0\ .$$
\end{lem}

\begin{proof} This is similar to section 25 of \cite{J2}. Denote by
$\Diff_f(M,A \rel R)$ the subgroup of $\Diff (M,A \rel R)$ consisting of
diffeomorphisms taking fibers to fibers.  The first assertion is that
the inclusion of $\Diff_f(M,A\rel R)$ into\break $ \Diff (M,A \rel R)$ induces
an isomorphism on $\pi_0$. A proof of this without the ``$\hbox{\it
rel}\,R$'' is indicated on pages~85--86 of \cite{W}, and the same proof
works $\hbox{\it rel\/}\, R$.

Since elements of $\Diff_f(M,A\rel R)$ take exceptional fibers to
exceptional fibers with the same labeling data, a natural homomorphism
$\Phi\co  \Diff_f(M,A \rel R) \to \Diff^*(S,E\cup p(A) \rel p(R))$ is
induced by projection. By theorem~8.3 of~\cite{Kalliongis-McCullough},
$\Phi$ is locally trivial, so is a Serre fibration. (One can also
check directly that $\Phi$ is a Serre fibration. Since we are dealing
with groups, it suffices to construct a $k$--parameter isotopy of the
identity of $M$ which lifts a given $k$--parameter isotopy of the
identity of $S$, and this is not difficult.) A section of $\Phi$ can
be constructed as follows.  Let $M_0$ be $M$ with an open fibered
tubular neighborhood of the exceptional fibers deleted, and let $S_0$
be the image of $M_0$ in $S$. If $S$ is orientable, then $M_0$ is a
product $S_0 \times S^1$, from which $M$ can be obtained by filling in
solid torus neighborhoods of the exceptional fibers in a standard way
depending only on the labeling data of the exceptional fibers.
Diffeomorphisms of $S_0$ give rise to diffeomorphisms of $M_0$ by
taking the product with the identity on $S^1$, and then these
diffeomorphisms extend over $M$ in the obvious way, assuming the
labeling data is preserved.  In case $S$ is non-orientable, $M_0$ can
be obtained by doubling the mapping cylinder of the orientable double
cover $\widetilde S_0 \to S_0$. Diffeomorphisms of $S_0$ lift
canonically to diffeomorphisms of $\widetilde S_0$, hence by taking
the induced diffeomorphisms of mapping cylinders we get a section of
$\Phi$ in this case too.

Thus from the exact sequence of homotopy groups of the fibration
$\Phi$ we obtain the split short exact sequence of the Proposition but
with $H_1(S,\bdy S -p(R))$ replaced by $\pi_0 (X)$ where $X$ is the
fiber of $\Phi$.  It remains then to produce an isomorphism $\pi_0 (X)
\cong H_1(S,\bdy S - p(R))$ (cf lemma~25.2 of~\cite{J2}).

The fiber $X$ of $\Phi$ consists of the diffeomorphisms taking each
circle fiber of $M$ to itself.  Note that orientations of fibers are
preserved since we are considering only diffeomorphisms which are the
identity on $R$. We may assume elements of $X$ restrict to rotations
of each circle fiber, in view of the fact that the groups of
orientation-preserving diffeomorphisms of $S^1$ has the homotopy type
of the rotation subgroup.  There is no harm in pretending the
exceptional fibers are not exceptional since the rotation of a
exceptional fiber is determined by the rotations of nearby fibers.  If
$S$ is orientable, then $M = S\times S^1$ and diffeomorphisms which
rotate fibers are the same as maps $(S,p(R)) \to (S^1,1)$, measuring
the angle of rotation in each fiber.  Thus $\pi_0 (X)$ is the group of
homotopy classes of maps $(S,p(R)) \to (S^1,1)$, ie $H^1(S,p(R))$,
which is isomorphic to $H_1(S,\bdy S -p(R))$ by duality.  When $S$ is
non-orientable one could presumably make the same sort of argument
using cohomology with local coefficients, but instead we give a direct
geometric argument, which applies when $S$ is orientable as well.

We can construct $S$ from a collar $p(R) \times I$ by attaching
1--handles, plus a single 2--handle if $p(R) = \bdy S$.  The core arcs of
the 1--handles, extended through the collar to $p(R)$, lift to annuli
$A_i$ in $M$ with $\bdy A_i \subset \partial M$.  Each diffeomorphism
in $X$ restricts to a loop of diffeomorphisms of $S^1$ on each $A_i$.
Since $\pi_1 (\Diff (S^1)) \cong \Z$, we thus have a homomorphism $\phi
\co\pi_0(X) \to \Z^n$ if there are $n \ A_i$'s.  Clearly $\phi$ is an
injection, so $\pi_0 X$ is finitely generated free abelian.  If $p(R) \ne
\bdy S$, so there is no 2--handle, then $\phi$ is obviously surjective
as well.  This is also true if $p(R) = \bdy S$, since it is not hard to
see that Dehn twist diffeomorphisms of the $A_i$'s extend to
diffeomorphisms in $X$.  A homology calculation shows that $H_1(S,\bdy
S-p(R))$ is a direct sum of $n$ copies of $\Z$ since $p(R)$ is non-empty.
\end{proof}

{\bf Remark}\stdspace The group $\hscr^*(S,E\cup p(A) \rel p(R))$ is
isomorphic to a subgroup of finite index in $\hscr(S_0,p(A) \rel
p(R))$, where $S_0$ is obtained from $S$ by deleting open discs about
the points of $E$. According to lemma~\ref{HARER}, the latter is
geometrically finite. Since $H_1(S,\partial S-p(R);\Z)$ is free
abelian when $R$ is non-empty, lemma~\ref{seifert mcg}, together with
lemma~\ref{higher homotopy} below, implies the Main Theorem in the
case when $M$ is a Seifert manifold.\medskip

The proof of our next lemma will use Dehn twists of $M$, which are
diffeomorphisms defined as follows. Let $F$ be a torus or annulus,
either properly imbedded in $M$ or contained in $\partial M$, and let
$F\times I$ be a submanifold of $M$ with $F=F\times\{0\}$. For a loop
$\gamma\co I\to \Diff(F)$ representing an element of
$\pi_1(\Diff(F),\hbox{\it id}_F)$, a Dehn twist is defined by putting
$h(y,t)=(\gamma_t(y),t)$ for $(y,t)\in F\times I$, and $h(x)=x$ for
$x\notin F\times I$.

\begin{lem}\sl\label{higher homotopy} Let $M$, $A$, and $R$ be as in
lemma~\ref{ibundle mcg} or~\ref{seifert mcg}. Then for all $i > 0$,
$\pi_i(\Diff(M,A\rel R))=0$.
\end{lem}

\begin{proof} Consider the restriction fibration
$$
\Diff(M \rel \bdy M) \longrightarrow \Diff(M,A\rel R) \buildrel \rho\over
\longrightarrow \Diff(\bdy M,A\rel R)\ .
$$
\noindent For $i\geq 1$, $\pi_i(\Diff(M \rel \bdy M))=0$ by \cite{H2}.
For $i>1$, $\pi_i(\Diff(\bdy M,A \rel R))=0$ by surface theory, so we
need only check injectivity of the boundary homomorphism
$\bdy\co\pi_1(\Diff(\bdy M,A\rel R))\to \pi_0(\Diff(M \rel \bdy M))
$. The only components of $\bdy M$ which can contribute to
$\pi_1(\Diff(\bdy M,A\rel R))$ are torus components disjoint from $R$.
Such a torus disjoint from $A$ contributes a $\Z\times\Z$ factor,
while a torus which contains components of $A$ contributes a $\Z$
factor. The boundary homomorphism takes these elements of
$\pi_1(\Diff(\bdy M,A\rel R))$ to Dehn twists supported near these
boundary tori. Since boundary tori are involved, we are in the
Seifert-fibered case, and we can assume these Dehn twists take fibers
to fibers. A non-zero element of the kernel of $\bdy$ would give a
non-trivial linear combination of these Dehn twists which was zero in
$\hscr(M\rel\bdy M)$. By projecting this linear combination onto
$\hscr(S \rel \bdy S) $ we see that it must be a linear combination of
Dehn twists taking each fiber to itself. But by our homology
interpretation of these Dehn twists, the only non-trivial combinations
which could be isotopically trivial are those involving twists near
all components of $\bdy M$. Since we are assuming $R$ is non-empty,
there are no such combinations in the image of~$\bdy$.
\end{proof}

\section{Hyperbolic Manifolds}
\label{hyperbolic}

Let $M$ be a compact connected orientable irreducible 3--manifold with
non-empty boundary. We decompose $\bdy M$ into three compact
subsurfaces meeting only in their boundary circles: $T$, the union of
the torus components of $\bdy M$; $A$, a disjoint union of annuli in
the other components; and $B$, the closure of $\bdy M-(A\cup T)$. We
assume that all the components of $B$ have negative Euler
characteristic. For brevity we write $C=A\cup T$, the ``cusps" of $M$.
Assume the following.
\begin{enumerate}
\item[(i)] $B$ and $C$ are incompressible in $M$.
\item[(ii)] Every $\pi_1$--injective map of a torus into $M$
is homotopic into~$T$.
\item[(iii)] Every $\pi_1$--injective map of pairs $(S^1\times
I,S^1\times\partial I)\to (M,B)$ is
homotopic as a map of pairs to a map carrying $S^1\times I$ into
either $A$ or~$B$.
\item[(iv)] $M$ is not homeomorphic to $S^1\times S^1\times I$.
\end{enumerate}
\noindent Note that assumptions (iii) and (iv) imply that $(M,A)$ is
not of the form $(F\times I,\partial F\times I)$.

Let $\Diff(M,A)$ denote the diffeomorphisms of $M$ that take $A$
diffeomorphically to $A$. These also must take $M-C$ to
$M-C$.

\begin{lem}\sl\label{HYPERBOLIC STRUCTURE} For $M$, $A$, $T$, and $C$
as above, $M-C$ has a hyperbolic structure with totally geodesic
boundary. Its group of isometries $\Isom(M-C)$ is finite, and
$\Isom(M-C)\to\hscr(M,A)$ is an isomorphism.
\end{lem}

\noindent The homomorphism $\Isom(M-C)\to\hscr(M,A)$ requires a bit of
explanation. Each component of $C$ inherits a Euclidean structure from
the corresponding cusp of $M-C$, and isometries of $M-C$ induce
isometries of these Euclidean annuli and tori. So each isometry of
$M-C$ extends uniquely to a diffeomorphism of $M$ preserving~$A$.

\begin{proof} Assume first that $T=\partial M$. By a celebrated result
of Thurston, $M-T$ has a complete hyperbolic structure of finite
volume, and by the Mostow Rigidity Theorem $\Isom(M-T)$ is finite and
the composition $\Isom(M-T)\to\hscr(M)\to \Out(\pi_1(M))$ is an
isomorphism. Since $M$ is aspherical, the outer automorphism group
$\Out(\pi_1(M))$ is naturally isomorphic to the group of homotopy
classes of homotopy equivalences from $M$ to $M$. Since every
incompressible torus in $M$ is homotopic into $\partial M$, an
application of the homotopy extension property shows that every
homotopy equivalence is homotopic to one which preserves $\partial M$.
By \cite{W} and the fact that $M$ is not of the form $F\times I$,
every boundary-preserving homotopy equivalence is homotopic to a
diffeomorphism. Therefore $\hscr(M)\to\Out(\pi_1(M))$ is surjective.
Also by \cite{W}, it is injective and the lemma follows in the case
$T=\partial M$.

Suppose now that $T\neq\partial M$. Let $N$ be the manifold obtained
by identifying two copies of $M$ along $B$ (using the identity map).
The boundary of $N$ is incompressible and consists of tori, and
assumption (iii) ensures that every incompressible torus in $N$ is
homotopic into $\partial N$. From the previous case, the interior of
$N$ admits a hyperbolic structure and $\Isom(N-\partial N)\to\hscr(N)$
is an isomorphism.  Let $\tau'$ be the involution of $N$ that
interchanges the two copies of $M$. Its fixed-point set is $B$.  Let
$\tau$ be the isometry in the isotopy class of $\tau'$. Note that
$\tau^2$ is an isometry isotopic to the identity, hence equals the
identity, so $\tau$ is an involution. By \cite{Tollefson}, homotopic
involutions of $N$ are strongly equivalent, ie there is a
homeomorphism $k$ of $N$, isotopic to the identity, such that $k\tau'
k^{-1}=\tau$. Regard $M$ as one of the copies of $M$ in $N$, so that
$B$ is its frontier. Then $k(B)$ is the fixed-point set of $\tau$, and
$k$ carries $M$ homeomorphically to the closure of one of the
complementary components of $k(B)$. Therefore by changing coordinates
using $k$ we may assume that $B$ is the fixed point set of $\tau$. The
fixed-point set of an isometry is totally geodesic so the restriction
to $M-C$ of the hyperbolic structure on $N-\partial N$ is complete
with totally geodesic boundary.

Define $\Phi\co \hscr(M,A)\to \hscr(N)$ by sending $\langle h\rangle$
to the class represented by $D(h)$, the double of $h$ along $B$. Let
${\cal T}$ be the subgroup of order~2 in $\hscr(N)$ generated by
$\langle\tau\rangle$. We claim that $\Phi$ is injective and that
$\Phi(\hscr(M,A))\times {\cal T}$ is the centralizer of
$\langle\tau\rangle$ in $\hscr(N)$. Suppose first that $\langle
h\rangle$ lies in the kernel of $\Phi$. Let $B_0$ be a component of
$B$. Assumption~(iii) implies that $N$ does not fiber over $S^1$ with
fiber $B_0$, so lemma~\ref{preserve F} implies that $D(h)$ is isotopic
to $1_N$ preserving $B_0$ at each level. Repeating for each component
of $B$, we find that $D(h)$ is isotopic to the identity preserving
$B$, so $\langle h\rangle$ was trivial in $\hscr(M, A)$. Clearly the
image of $\Phi$ lies in the centralizer, since $D(h)$ actually
commutes with $\tau$. Suppose $\langle H\rangle$ is an element in the
centralizer of $\tau$. Then $H\tau H^{-1}$ is isotopic to $\tau$.
Again by Tollefson's result, they must be strongly equivalent. So $H$
is isotopic to $kH$ with $kH \tau (kH)^{-1}=\tau$. This implies that
$kH$ preserves $B$. If $kH$ does not reverse the sides of $B$, then it
must be of the form $D(h)$, so lies in $\Phi(\hscr(M,A))$. It it does
reverse the sides, then $kH\tau$ does not. If follows that
$\Phi(\hscr(M,A))$ and ${\cal T}$ generate the centralizer. If $kH$
and $k'H$ preserve $B$ and are isotopic, then by lemma~\ref{preserve
F} they are isotopic preserving $B$. Therefore it is well-defined
whether an element in the centralizer of $\langle\tau\rangle$ in
$\hscr(N)$ preserves the sides of $B$. In particular, elements of the
image of $\Phi$ do not reverse the sides of $B$, so
$\Phi(\hscr(M,A))\cap {\cal T}$ consists only of the identity, and the
claim follows.

From the case $T=\partial M$, $\Isom(N-\partial N)\to\hscr(N)$ is an
isomorphism. An isometry on $N$ commutes with $\tau$ and preserves the
sides of $B$ if and only if it is the double along $B$ of an isometry
of $M$. Therefore sending an isotopy class in $\Phi(\hscr(M,A))$ to
the restriction to $M$ of the unique isometry that it contains defines
an inverse to $\Isom(M-C)\to\hscr(M,A)$.
\end{proof}

\begin{prop}\sl\label{boundary twists} Let $M$, $B$, $A$, and $T$ be
as above. Let $R$ be a non-empty union of components of $B$, $A$,
and $T$, and let $R_0$ be the components of $R$ that have Euler
characteristic zero. Then $\hscr(M,A\rel R)\cong H_1(R_0;\Z)$, and
$\pi_i(\Diff(M,A\allowbreak\rel R))=0$ for $i\geq 1$.
\end{prop}

\begin{proof} Consider the fibration $\Diff(M,A\rel R)\to \Diff_R(M,A)
\allowbreak\to \diff(R)$ where $\diff(R)$ is the identity component of
$\Diff(R)$ and $\Diff_R(M,A)$ is the subgroup of $\Diff(M,A)$
consisting of diffeomorphisms taking each component of $R$ to itself
by a diffeomorphism isotopic to the identity.  This fibration gives an
exact sequence:
$$ \pi_1 \diff(R) \buildrel \bdy \over \longrightarrow  \hscr
(M,A\,rel\,R) \longrightarrow
\hscr_R(M,A) \longrightarrow0 \leqno (*)
$$

\noindent
The following argument shows that the map $\bdy$ is injective, so
$(*)$ is in fact a short exact sequence. First note that $\diff(R)$ is
the direct product of the $\diff(F)$ as $F$ ranges over the components
of $R$. For the components that have negative Euler characteristic,
$\diff(F)$ is contractible, while if $F$ is a torus or annulus,
$\diff(F)$ is homotopy equivalent to $F$. If $R_0$ is empty, then
$\pi_1(\diff(R))$ is trivial. Otherwise, fix a component $F$ of $R_0$.
The boundary map $\bdy$ sends elements $\langle \gamma\rangle$ of
$\pi_1(\diff(F))$ to Dehn twists supported in a collar neighborhood of
$F$. Such a Dehn twist induces an inner automorphism of $\pi_1(M,x_0)$
for $x_0 \in F$, namely, conjugation by the element of $\pi_1(M,x_0)$
represented by the loop in $F$ around which a basepoint $x_0$ in $F$
is carried by $\gamma$. This element uniquely determined
$\langle\gamma\rangle$, in particular it is non-trivial when $\gamma$
is non-trivial. Note that Dehn twists near other components of $R_0$
have no effect on $\pi_1(M,x_0)$. Inner automorphisms of $\pi_1 M$ are
always non-trivial since $\pi_1 M$ has trivial center (a standard fact
about hyperbolic 3--manifolds other than the ones ruled out by
assumptions~(i)--(iv) above). Therefore the map $\pi_1 \diff(F) \to
\hbox{\it Aut}(\pi_1(M,x_0))$ is injective. Fixing basepoints
$x_1,\ldots\,$, $x_k$ in the components of $R_0$, we obtain a
composition $\pi_1 \diff(R) \buildrel \bdy \over \rightarrow  \hscr
(M\,rel\,R)\rightarrow \prod_{i=1}^k\hbox{\it Aut}(\pi_1(M,x_i))$
which is injective, showing that $\bdy$ is injective.

Next we show that $\pi_i(\Diff(M,A\rel R))=0$ for $i\geq1$. It is
sufficient as in the proof of lemma~\ref{higher homotopy} to check
that $\partial\co \pi_1(\Diff(\bdy M,A\rel R))\to{\cal H}(M\rel
\partial M)$ is injective. Since $\pi_1(\diff(F))$ is trivial when $F$
has negative Euler characteristic, $\pi_1(\Diff(\bdy M,A\rel R))$ is
generated by elements of $\pi_1(\diff(F))$ for the torus components of
$\partial M$ that are not contained in $R$. Fixing basepoints $y_j$ in
these components, we have as before an injective homomorphism
$$\pi_1(\Diff(\partial M,A\rel R) \buildrel \bdy \over \longrightarrow
{\cal H}(M\rel \partial M)\to\prod\hbox{\it Aut}(\pi_1(M,y_j)),$$
showing that $\partial$ is injective.

Now we turn to the calculation of $\hscr(M,A\rel R)$. By
lemma~\ref{HYPERBOLIC STRUCTURE}, we can fix a hyperbolic structure on
$M-(A\cup T)$ such that $\Isom(M-(A\cup T))$ is finite and
$\Isom(M-(A\cup T))\to \hscr(M,A)$ is an isomorphism. Suppose first
that $R\neq R_0$. Since no non-trivial isometry of a hyperbolic surface
of negative Euler characteristic is isotopic to the identity, the
subgroup of $\Isom(M-(A\cup T))$ that maps to the subgroup
$\hscr_R(M,A)$ of $\hscr (M,A)$ is trivial if $R\neq R_0$, in which
case $(*)$ gives ${\cal H}(M,A\rel R)\cong \pi_1(\diff(R))\cong
H_1(R_0;\Z)$.

Thus we may assume that $R=R_0$. The subgroup $\Isom\,_R(M-(A\cup T))
\subset\Isom(M-(A\cup T))$ that corresponds to the subgroup
$\hscr_R(M,A)\subset\hscr (M,A)$ consists of isometries which on each
component of $R$ are rotations isotopic to the identity. For each $
\varphi_0 \in \Diff(M,A \rel R) $ there is an isotopy $ \varphi_t$ in
$\Diff(M,A) $ from $\varphi_0$ to the isometry
$\varphi_1\in\Isom\,_R(M-(A\cup T))$ corresponding to $\varphi_0$
under the map $\hscr(M,A\rel R) \to \hscr_R(M,A)$. The isotopy
$\varphi_t$ is unique up to deformation since the fact that
$\pi_1\Diff(M,A\rel R)=0$ implies that  $\pi_1\Diff(M,A)=0$ by looking
a few terms to the left in the sequence $(*)$.

The group of rotations of $R$ can be identified with a subspace $R^*
\subset R$ by evaluation of rotations at a chosen basepoint in each
component of $R$; in annulus components we choose the basepoint in the
boundary of the annulus. The inclusion $R^*\hookrightarrow R$ is a
homotopy equivalence. Let $G$ be the subgroup of $R^*$ obtained by
restriction of $\Isom\,_R(M-(A\cup T))$. Evaluation of the path
$\varphi_t$ at the basepoints in $R$ then gives a well-defined map
$\Phi\co\hscr(M,A\rel R)\to \pi_1(R^*,G)$ which is a homomorphism
if $\pi_1(R^*,G)$ is given the group structure induced by the group
structure of $R^*$. In fact, $\Phi$ gives a map from the short exact
sequence $(*)$ to the short exact sequence
$0\to\pi_1(R)\to\pi_1(R^*,G)\to G\to 0$, hence $\Phi$ is an
isomorphism by the five-lemma. By lifting paths to the universal cover
of $R^*$, we identify $\pi_1(R^*,G)$ with a cocompact lattice in a
Euclidean space, containing the deck transformation group $\pi_1(R)$
as a subgroup of finite index. Thus the group $\hscr(M,A \rel R)\cong
\pi_1(R^*,G)$ is abstractly isomorphic to $\pi_1(R) \cong H_1(R;\Z)$.
\end{proof}

{\bf Remark}. The case when $A$ is empty yields the Main Theorem in
the case when $M$ is a simple manifold.

\section{Decomposable Manifolds}
\label{decomposable}

In this section we prove the Filtration Theorem in the general case.
The Main Theorem follows immediately using lemmas~\ref{extensions of
gf groups} and~\ref{HARER}.

Suppose first that $\partial M$ is compressible. We assume the
compressible components of $\bdy M$ lie in $R$. By inductive
application of the Loop Theorem, one can construct a collection $E$ of
finitely many disjoint properly-imbedded discs, none of which is
isotopic into $\partial M$, such that each component of $M$ cut along
$E$ has incompressible boundary. Note that $\pi_0\Imb(E,M\rel \partial
E)=0$. By lemma~\ref{preserve F}, $\pi_1(\Imb(E,M\rel \partial E))=0$,
so the restriction fibration
$$
\Diff(M\rel E\cup R)\to \Diff(M\rel R)\to\Imb(E,M\rel\partial E)
$$
\noindent shows that $\hscr(M\rel R)\cong\hscr(M\rel E\cup R)$. The
latter group can be identified with $\hscr(M'\rel R')$, where $M'$ is
the result of cutting $M$ along $E$ and $R'$ is the union of boundary
components of $M'$ corresponding to $R$. Although $M'$ may no longer
be connected, each of its components meets $R'$. Thus we reduce to the
case that $R$ is incompressible. In particular, if $M$ was a
handlebody, then $\hscr(M\rel\partial M)$ is trivial and the proof is
completed.

Assuming now that $M$ has incompressible boundary, the elementary form
of the Torus--Annulus Decomposition Theorem of Jaco--Shalen and
Johannson (\cite{JS,J2}) says $M$ contains a 2--dimensional submanifold
$T\cup A$, where $T$ consists of incompressible tori and $A$ of
incompressible annuli, such that each component $W$ of the manifold
obtained by splitting $M$ along $T\cup A$ is either:

\begin{enumerate}
\item[(a)] simple, meaning that $W$, $T_W$, and $A_W$ satisfy
conditions (i)--(iv) at the beginning of section~\ref{hyperbolic},
where $T_W$ is the union of the torus boundary components of $W$ and
$A_W$ is the union of the components of the closure of ${\bdy W - \bdy
M}$ that are annuli;
\item[(b)]  an $I$--bundle over a surface of negative Euler
characteristic, such that
$W \cap \bdy M$ is the associated $\partial I$--bundle; or
\item[(c)]  Seifert-fibered, with $ W \cap \bdy M$  a union of fibers.
\end{enumerate}

\noindent
Further, when $T\cup A$ is chosen to be minimal with respect to
inclusion among all such submanifolds, it is unique up to ambient
isotopy of $M$. We need a relative form of this uniqueness: Two
choices of $T\cup A$ having the same boundary are isotopic fixing $R$.
This follows from the previous uniqueness statement since the
obstruction to deforming an arbitrary isotopy to an isotopy fixing $R$
is the homotopy class of the loop of embeddings of $R\cap\bdy A$
traced out during the isotopy, but  $\bdy A$ is disjoint from torus
components of $\bdy M$ (an exercise from the definitions) and there
are no non-trivial loops of embedded circles in surfaces of negative
Euler characteristic.

This relative uniqueness implies that the natural map $\hscr(M,T\cup A
\rel R) \to \hscr(M\rel R)$ is surjective. By lemma~\ref{preserve F}
it is also injective, so we have $\hscr(M,T\cup A\rel R) \cong
\hscr(M\rel R)$. This remains true if we replace each of the annuli
and tori of $T\cup A$ by two nearby parallel copies of itself, still
calling the doubled collection $T\cup A$. The advantage in doing this
is that now when we split $M$ along $T\cup A$, the pieces produced by
the splitting are  submanifolds of $M$, and $M$ is their union. The
new pieces lying between parallel annuli and tori of the original
$T\cup A$ we view as additional Seifert-fibered pieces.

Let $V_1=W_1$ be a piece which meets $R$, and inductively, let $V_i =
V_{i-1}\cup W_i$ where $W_i$ is a piece which meets $V_{i-1}$, other
than the pieces already in $V_{i-1}$. For completeness let $V_0$ be
empty. Then we have restriction fibrations
$$
\Diff(M,T\cup A\rel V_i\cup R)\longrightarrow
\Diff(M,T\cup A\rel V_{i-1}\cup R)\buildrel \rho
\over \longrightarrow\Diff(W_i,A_i\rel R_i)
$$
where $A_i=W_i\cap A$ and $R_i=W_i\cap(V_{i-1}\cup R)$. These fibrations
yield exact sequences
$$
0\to\hscr(M,T\cup A\rel V_i\cup R)\rightarrow
\hscr(M,T\cup A\rel V_{i-1}\cup R)\buildrel \rho_*
\over \longrightarrow\hscr(W_i,A_i\rel R_i)
$$
where the zero at the left end is $\pi_1(\Diff(W_i,A_i\rel R_i))$,
which vanishes by lemma~\ref{higher homotopy} or
proposition~\ref{boundary twists}. The Filtration Theorem will follow
once we show that the image of each map $\rho_*$ has a filtration of
the sort in the theorem.

The case that $W_i$ is hyperbolic is immediate since
$\hscr(W_i,A_i\rel R_i)$ is finitely generated free abelian by
proposition~\ref{boundary twists}, hence also any subgroup of it.
Consider next the case that $W_i$ is an $I$--bundle. By
lemma~\ref{ibundle mcg}, $\hscr(W_i,A_i\rel R_i)=0$ unless $R_i\subset
A_i$. In the latter case $\hscr(W_i,A_i\rel R_i)\cong\hscr(S\rel
p(R))$ with $p(R)$ a union of components of $\partial S$. The image of
$\rho_*$ has finite index in this group since it contains the
finite-index subgroup represented by diffeomorphisms of $W_i$ which
restrict to the identity on $A_i$, corresponding to elements of
$\hscr(S\rel p(R))$ represented by diffeomorphisms which are the
identity on $\bdy S$.

There remains the case that $W_i$ is Seifert-fibered. By
lemma~\ref{seifert mcg}, there is a short exact sequence
$$\displaylines{
0   \longrightarrow   H_1 (S,\bdy S -p(R_i);\Z)   \longrightarrow
\hscr(W_i,A_i\rel R_i)    \longrightarrow \cr
\hscr^*(S,E\cup p(A_i) \rel p(R_i))\longrightarrow 0\ .\cr}
$$
As we noted in the Remark at the end of section 2, $\hscr^*(S,E\cup
p(A_i) \rel p(R_i))$ is a subgroup of finite index in a mapping class
group ${\cal H}(S_0,p(A_i)\rel p(R_i))$. The image of $\rho_*$
projects into $\hscr^*(S,E\cup p(A_i) \rel p(R_i))$ as a subgroup of
finite index, since the image of $\rho_*$ contains the isotopy classes
represented by diffeomorphisms which are the identity on $\bdy W_i$.
Since $R_i$ is not empty, $H_1(S,\partial S-p(R_i);\Z)$ is free
abelian. So it intersects the image of $\rho_*$ in a free abelian
group, and the proof is complete.

%
%

%
%
\end{document}